\RequirePackage{fix-cm}
\documentclass[smallextended]{svjour3}       % onecolumn (second format)
\smartqed  % flush right qed marks, e.g. at end of proof
\usepackage{graphicx}
\usepackage{amssymb}
\usepackage{amsmath}
\usepackage[ruled, lined]{algorithm2e}
\usepackage{subfig}
\begin{document}

\title{
A Fast Implementation of Singular Value Thresholding Algorithm using Recycling Rank Revealing Randomized Singular Value Decomposition%\thanks{Grants or other notes
}
% \subtitle{Do you have a subtitle?\\ If so, write it here}

\titlerunning{Fast Implementation of Singular Value Thresholding}        % if too long for running head

\author{Yaohang Li         \and
        Wenjian Yu %etc.
}


\institute{Yaohang Li (Corresponding Author)\at
              Department of Computer Science, Old Dominion University, Norfolk, VA 23529, USA\\
              \email{yaohang@cs.odu.edu}           %  \\
           \and
           Wenjian Yu \at
              Department of Computer Science and Technology, Tsinghua National Lab of Information Science and Technology, Tsinghua University, Beijing 100084, China\\
              \email{yu-wj@tsinghua.edu.cn}           %  \\ 
}

\date{Received: date / Accepted: date}

\maketitle

\begin{abstract}
In this paper, we present a fast implementation of the Singular Value Thresholding (SVT) algorithm for matrix completion. A rank-revealing randomized singular value decomposition (R\textsuperscript{3}SVD) algorithm is used to adaptively carry out partial singular value decomposition (SVD) to fast approximate the SVT operator given a desired, fixed precision. We extend the R\textsuperscript{3}SVD algorithm to a recycling rank revealing randomized singular value decomposition (R\textsuperscript{4}SVD) algorithm by reusing the left singular vectors obtained from the previous iteration as the approximate basis in the current iteration, where the computational cost for partial SVD at each SVT iteration is significantly reduced. A simulated annealing style cooling mechanism is employed to adaptively adjust the low-rank approximation precision threshold as SVT progresses. Our fast SVT implementation is effective in both large and small matrices, which is demonstrated in matrix completion applications including image recovery and movie recommendation system. 
\end{abstract}

\section{Introduction}
\label{intro}
Given a matrix $\mathbf{A}\in\mathbb{R}^{m \times n}$ with missing entries, the objective of matrix completion \cite{Ref18} is to recover the matrix by inferring the missing ones from the known set of entries $(i,j)\in \mathbf{\Lambda}$. The formulation of the mathematical model for matrix completion is low-rank assumption, i.e., $\mathbf{A}$ has rank ${r \ll m, n}$. Under the low-rank assumption, there exist an $m \times r$ matrix $\mathbf{M}$ and an $r \times n$ matrix $\mathbf{N}$ such that $\mathbf{MN}=\mathbf{A}$. Then, the matrix completion problem becomes an optimization problem such that 
\begin{equation}
\begin{split}
&\min_{\mathbf{X}} rank(\mathbf{X}) \\
&s.t., \mathbf{X}_{ij}=\mathbf{A}_{ij}, (i,j) \in \mathbf{\Lambda}.
\end{split}
\end{equation}
Unfortunately, the above rank minimization problem is known to be NP-hard and thus is impractical for large-scale matrix completion problems. A relaxation form is proposed by minimizing the sum of the singular values of $\mathbf{A}$, which is known as the nuclear norm of $\mathbf{A}$. Then, the matrix completion problem is reformulated as a convex optimization problem \cite{Ref17} such as
\begin{equation}
\begin{split}
&\min_{\mathbf{X}} \|\mathbf{X}\|_* \\
&s.t., \mathbf{X}_{ij}=\mathbf{A}_{ij}, (i,j) \in \mathbf{\Lambda},
\end{split}
\end{equation}
where $\|.\|_*$ denotes the nuclear norm. Candes and Recht \cite{Ref2} have shown that under certain conditions, the solution obtained by optimizing the nuclear norm is equivalent to the one by rank minimization. 

The nuclear norm optimization problem for matrix completion can be efficiently addressed by using the singular value thresholding (SVT) algorithm \cite{Ref1}, which is a first-order algorithm approximating the nuclear norm optimization problem by
\begin{equation}
\begin{split}
&\min_{\mathbf{X}} {1 \over 2}\|\mathbf{X}\|_F^2 + \tau \|\mathbf{X}\|_* \\
&s.t., \mathbf{X}_{ij}=\mathbf{A}_{ij}, (i,j) \in \mathbf{\Lambda} 
\end{split}
\end{equation}
with a threshold parameter $\tau$. Then, starting from an initial matrix $\mathbf{Y}^{(0)}$, where $\mathbf{Y}_{ij}^{(0)}=\mathbf{A}_{ij}$ for $(i,j) \in \mathbf{\Lambda}$ and $\mathbf{Y}_{ij}^{(0)}=0$ for $(i,j)\notin \mathbf{\Lambda}$, SVT applies an iterative gradient ascend approach formulated as Uzawa’s algorithm \cite{Ref9} or linearized Bregman’s iterations \cite{Ref10} such that

\begin{equation}
\label{Bregman}
\begin{split}
&\mathbf{X}^{(i)}=D_{\tau}(\mathbf{Y}^{(i)}) \\
&\mathbf{Y}^{(i+1)}=\mathbf{Y}^{(i)}+\delta P_{\mathbf{\Lambda}}(\mathbf{A}-\mathbf{X}^{(i)}),
\end{split}
\end{equation}
where $\delta$ is the step size, $P_\mathbf{\Lambda}$ is an orthogonal projector onto $\mathbf{\Lambda}$, and $D_\tau$ is known as the SVT operator. Given $\mathbf{Y}^{(i)}$ at the $i$th SVT iteration step and its singular value decomposition (SVD) $\mathbf{Y}^{(i)}=\mathbf{U}^{(i)} \mathbf{\Sigma}^{(i)} {\mathbf{V}^{(i)}}^T$, where $\mathbf{U}^{(i)}$ and $\mathbf{V}^{(i)}$ are orthonormal matrices and $\mathbf{\Sigma}^{(i)}=diag(\sigma_1^{(i)},\sigma_2^{(i)},…,\sigma_r^{(i)})$ is a diagonal matrix with $\sigma_1^{(i)} \geq \sigma_2^{(i)} \geq ,…, \geq \sigma_r^{(i)} \geq 0$ as the singular values of $\mathbf{Y}^{(i)}$, the SVT operator $D_{\tau}(\mathbf{Y}^{(i)})$ is defined as shrinking the singular values less than $\tau$ as well as their associated singular vectors, i.e.,
\begin{equation}
D_{\tau}(\mathbf{Y}^{(i)})= \sum_{j}^{\sigma_j^{(i)}\geq \tau}(\sigma_j^{(i)}-\tau)u_j^{(i)}{v_j^{(i)}}^T, \\
\end{equation}
where $u_j^{(i)}$ and $v_j^{(i)}$ are column vectors in $\mathbf{U}^{(i)}$ and $\mathbf{V}^{(i)}$, respectively.

Computing $D_{\tau}(\mathbf{Y}^{(i)})$ is the main operation in SVT, which is required to be repeatedly carried out at every iteration. A straightforward way to estimate $D_{\tau}(\mathbf{Y}^{(i)})$ is to compute full SVD on $\mathbf{Y}^{(i)}$ and then shrink the small singular values below threshold. However, for a relatively large matrix $\mathbf{Y}^{(i)}$, computing full SVD is costly, which prevents SVT from scaling up to large matrix completion problems. There are two different strategies to reduce the computational cost of evaluating the SVT operators. One is to replace full SVD with matrix operations of lower computational cost. For example, Cai and Osher \cite{Ref3} reformulate $D_{\tau}(\mathbf{Y}^{(i)})$ by projecting $\mathbf{Y}^{(i)}$ onto a 2-norm ball and then apply complete orthogonal decomposition and polar decomposition to the projection to obtain $D_{\tau}(\mathbf{Y}^{(i)})$, which saves 50\% or more computational time compared to SVT using full SVD. A more popular strategy is to compute partial SVD using Krylov subspace algorithms for the singular values of interest. This is due to the fact that only those singular values exceeding $\tau$ and their associated singular vectors in $\mathbf{Y}^{(i)}$ are concerned at each SVT iteration step \cite{Ref15}\cite{Ref16}. As a result, Krylov subspace algorithms can efficiently compute partial SVD that only reveals the singular values of interest. In fact, several SVT implementations \cite{Ref1}\cite{Ref15}\cite{Ref19}\cite{Ref20} use Lanczos algorithm with partial re-orthogonalization provided by PROPACK \cite{Ref4}. Combined with a rank prediction mechanism that can predict the number of singular values/vectors needed, partial SVD can efficiently accelerate SVT operator calculation if the number of singular values over $\tau$ is small compared to $\min(m,n)$. Nevertheless, the computational cost of Krylov subspace partial SVD relies on the number of singular values/vectors that need to be computed. As shown in \cite{Ref5}\cite{Ref26}, if the number of singular values/vectors needed exceeds $0.2\min(m,n)$, computing partial SVD is usually even more costly than carrying out full SVD.

In this paper, we design randomized algorithms to adaptively carry out partial SVD in order to fast approximate SVT operator computation. The fundamental idea of the randomized SVD algorithms \cite{Ref8} is to condense the sparse matrix $\mathbf{Y}^{(i)}$ into a small, dense matrix by projecting $\mathbf{Y}^{(i)}$ onto a sampling matrix as an approximate basis while keeping the important information of $\mathbf{Y}^{(i)}$. Then, performing a deterministic SVD on this small, dense matrix can be used to approximate the top singular values/vectors of $\mathbf{Y}^{(i)}$. Here, we first use a rank-revealing randomized singular value decomposition (R\textsuperscript{3}SVD) algorithm to adaptively carry out partial SVD to fast approximate $D_{\tau}(\mathbf{Y}^{(i)})$ under a fixed precision. Moreover, the accuracy of the randomized SVD algorithm relies on the quality of the sampling matrix - if the projection of $\mathbf{Y}^{(i)}$ onto the approximate basis captures a significant portion of the actions in $\mathbf{Y}^{(i)}$, it can lead to a good approximation of the top singular values/vectors of $\mathbf{Y}^{(i)}$ that can be used to compute the SVT operator. Therefore, based on R\textsuperscript{3}SVD, we design a recycling rank-revealing randomized SVD (R\textsuperscript{4}SVD) algorithm according to the fact that during SVT iterations, the change between the subsequent matrices $\mathbf{Y}^{(i)}$ and $\mathbf{Y}^{(i+1)}$ is relatively small, particularly at the later phase. As a result, the singular vectors obtained from the previous iteration step can potentially be recycled to build up the approximate basis in the next iteration during matrix sampling so that the computational cost of the subsequent SVT operator can be significantly reduced. In R\textsuperscript{4}SVD, starting from the left singular vectors of $\mathbf{Y}^{(i)}$ obtained from the previous iteration step, the sampling matrix is incrementally built up while the approximation error of $\mathbf{Y}^{(i+1)}$ projected onto the approximate basis derived from the sampling matrix is monitored to determine the appropriate rank to satisfy the given precision. We also adopt a simulated annealing-style cooling scheme \cite{Ref13}\cite{Ref14} to adaptively adjust the required partial SVD precision along SVT iterations. The computational efficiency of our fast SVT implementation is demonstrated on the applications of image recovery and movie recommendation system.

\section{Description of Algorithms}
\label{sec:2}
The fundamental idea of our fast SVT implementation is to use randomized SVD algorithms to fast approximate partial SVD in order to accelerate the SVT operator $D_{\tau}(\mathbf{Y}^{(i)})$ at the $i$th SVT iteration. This is a fixed precision problem, i.e., given the desired error percentage threshold $\epsilon$, the goal is to minimize the rank parameter $k$ of a low-rank SVD approximation $\mathbf{Y}_L^{(i)}=\mathbf{U}_L^{(i)} \mathbf{\Sigma}_L^{(i)} {\mathbf{V}_L^{(i)}}^T$, so that $\|\mathbf{Y}_L^{(i)}-\mathbf{Y}^{(i)}\|/\|\mathbf{Y}^{(i)}\| \leq \epsilon_{threshold}$. Nevertheless, the original randomized SVD (RSVD) algorithm is designed to address a fixed rank problem, i.e., given the fixed rank parameter $k$, the goal of RSVD is to obtain a $k$-rank SVD approximation $\mathbf{Y}_L^{(i)}=\mathbf{U}_L^{(i)} \mathbf{\Sigma}_L^{(i)} {\mathbf{V}_L^{(i)}}^T$ in order to minimize $\|\mathbf{Y}_L^{(i)}-\mathbf{Y}^{(i)}\|/\|\mathbf{Y}^{(i)}\|$. To address the fixed precision problem, we first present a rank revealing randomized SVD (R\textsuperscript{3}SVD) algorithm where the appropriate rank parameter to satisfy the fixed precision is gradually revealed by incrementally building up the low-rank SVD approximation while estimating the approximation error. Then, we introduce a Recycling Rank-Revealing Randomized SVD (R\textsuperscript{4}SVD) algorithm that can take advantage of the singular vectors obtained from previous iterations to accelerate the randomized SVD computation. Finally, by putting all pieces together, we present our fast implementation of SVT algorithm.

\subsection{Randomized SVD (RSVD)}
\label{sec:2.1}
The RSVD method was proposed by Halko et al. \cite{Ref8} to approximate the top $k$ singular values and singular vectors of a given matrix $\mathbf{Y}^{(i)}$, where $k$ is a given fixed rank. The fundamental idea of the RSVD algorithm using Gaussian random sampling is to construct a small condensed subspace from the original matrix, where the dominant actions of $\mathbf{Y}^{(i)}$ could be quickly estimated from this small subspace with relatively low computation cost and high confidence. The procedure of RSVD \cite{Ref8} is described as follows.

\begin{algorithm}
 \KwIn{$\mathbf{Y}^{(i)} \in \mathbb{R}^{m \times n}$, number of power iteration $np \in \mathbb{N}$, $k \in \mathbb{N}$ and $p \in \mathbb{N}$ satisfying $k+p \leq \min(m,n)$}
 \KwOut{$\mathbf{U}_L^{(i)} \in \mathbb{R}^{m \times k}$, $\mathbf{\Sigma}_L^{(i)} \in \mathbb{R}^{k \times k}$, and $\mathbf{V}_L^{(i)} \in \mathbb{R}^{k \times n}$ }
 \bigbreak
 Construct an $n \times (k+p)$ Gaussian random matrix $\mathbf{\Omega}$\;
 $\mathbf{X} \gets \mathbf{Y}^{(i)}\mathbf{\Omega}$\;
 \For {$j\leftarrow 1$ \KwTo $np$} { 
 $\mathbf{X} \gets \mathbf{Y}^{(i)}\mathbf{X}$ \tcc*{power iterations}
 }
 $[\mathbf{Q},\mathbf{R}] \gets qr(\mathbf{X})$ \tcc*{construct orthogonal subspace}
 $\mathbf{B} \gets \mathbf{Q}^T \mathbf{Y}^{(i)}$ \tcc*{QB decomposition}  
 $[\mathbf{U}_B, \mathbf{\Sigma}_B, \mathbf{V}_B] \gets svd(\mathbf{B})$\;
 $\mathbf{U}_B \gets \mathbf{Q}\mathbf{U}_B$\;
 $\mathbf{U}_L^{(i)} \gets \mathbf{U}_B(:,1:k)$, $\mathbf{\Sigma}_L^{(i)} \gets \mathbf{\Sigma}_B(1:k,1:k)$, and $\mathbf{V}_L^{(i)} \gets \mathbf{V}_B(:,1:k)$
 \caption{Randomized SVD Algorithm with Gaussian Sampling}
\end{algorithm}

Starting from a Gaussian random matrix $\mathbf{\Omega}$, the RSVD algorithm projects the original matrix $\mathbf{Y}^{(i)}$ onto $\mathbf{\Omega}$ as $\mathbf{X} \leftarrow \mathbf{Y}^{(i)}\mathbf{\Omega}$. Let $\omega_j\in\mathbb{R}^n$ and $x_j \in \mathbb{R}^n$ denote the $j$th column vector of random matrix $\mathbf{\Omega}$ and the $j$th column vector of matrix $\mathbf{X}$, respectively. Since each element in $\mathbf{\Omega}$ is chosen independently, $\omega_j$ can be represented as
\begin{equation}
\omega_j=c_{1j}v_1+c_{2j}v_2+...+c_{nj}v_n, \textrm{for }j=1,...,k+p\\
\end{equation}
where $v_l \in \mathbb{R}^n$ is $l$th right singular vector of matrix $\mathbf{Y}^{(i)}$ and $c_{lj} \neq 0$ with probability $1.0$. In RSVD, after simply projecting $\mathbf{Y}^{(i)}$ onto $\mathbf{\Omega}$, we could have
\begin{equation}
x_j=\sigma_1 c_{1j} v_1 + \sigma_2 c_{2j} v_2 + ... + \sigma_n c_{nj} v_n, 
\end{equation}
where $\sigma_l$ is the $l$th singular value of $\mathbf{Y}^{(i)}$ sorted by non-decreasing order such that $\sigma_1^{(i)} \geq \sigma_2^{(i)} \geq ,..., \geq \sigma_t^{(i)} \geq 0$ and $\sigma_l c_{lj}$ constitutes the weight of $x_j$ on $v_i$. Consequently, Gaussian sampling ensures that all singular vectors are kept in the subspace but the singular vectors corresponding to larger singular values likely yield bigger weights in $x_j$. Therefore, compared to $\omega_j$, weights of the dominant right singular vectors are amplified by the corresponding singular values. As a result, the space spanned by the columns of $\mathbf{X}$ reflects dominating weights in high probability on the singular vectors corresponding to the top $k$ singular values. Moreover, for stability consideration, an oversampling parameter $p$ is used to serve as a noise-filter to get rid of unwanted subspace corresponding to relative small singular values, when the SVD decomposition on $\mathbf{B}$ is carried out to approximate the top $k$ singular values/vectors of $\mathbf{Y}^{(i)}$. In practice, $p$ is given with a small value, such as $5$ or $10$, as suggested by Halko et al. \cite{Ref8}.

Compared to full SVD directly operating on the $m \times n$ matrix $\mathbf{Y}^{(i)}$, which is rather computational costly when both $m$ and $n$ are large, the major operations in RSVD are carried out
on the block matrices instead. These block matrix operations include matrix-block matrix multiplications as well as QR and SVD on the block matrices. Specifically, matrix-block matrix multiplications between $\mathbf{Y}^{(i)}$ and the block matrices take $O((np+2)(k+p)T_{mult})$ floating-point operations, where $T_{mult}$ denotes the computational cost of a matrix-vector multiplication. 

Notice that given the fixed rank parameter $k$, the goal of RSVD is to obtain a $k$-rank SVD approximation $\mathbf{Y}_L^{(i)}=\mathbf{U}_L^{(i)} \mathbf{\Sigma}_L^{(i)} {\mathbf{V}_L^{(i)}}^T$ to minimize $\|\mathbf{Y}_L^{(i)}-\mathbf{Y}^{(i)}\|/\|\mathbf{Y}^{(i)}\|$, which is a fixed rank problem. In SVT, the appropriate rank parameter $k$ to obtain a partial SVD approximation with desired precision is unknown beforehand. This is a fixed precision problem. The R\textsuperscript{3}SVD algorithm described in Section \ref{sec:2.2} is designed to address the fixed precision problem by adaptively determining the approximate rank parameter $k$ and incrementally build up the low-rank SVD approximation.
\subsection{Rank-Revealing Randomized SVD (R\textsuperscript{3}SVD)}
\label{sec:2.2}
The R\textsuperscript{3}SVD algorithm is designed to address the fixed precision problem. Compared to RSVD, R\textsuperscript{3}SVD incorporates three major changes including orthogonal Gaussian sampling, adaptive QB decomposition, and stopping criteria based on error percentage estimation. Moreover, since the sampling error is precisely estimated, oversampling in the RSVD is no longer necessary in R\textsuperscript{3}SVD. To illustrate the R\textsuperscript{3}SVD algorithm, for a given matrix $\mathbf{Y}^{(i)} \in \mathbb{R}^{m \times n}$ at the $i$th iteration of SVT and its $k$-rank approximation $\mathbf{Y}_L^{(i)}$, we first define the error percentage of $\mathbf{Y}_L^{(i)}$ with respect to $\mathbf{Y}^{(i)}$ measured by the square of Frobenius norm, i.e.,
\begin{equation}
\epsilon={{{\|\mathbf{Y}_L^{(i)}-\mathbf{Y}^{(i)}\|}_F^2} \over {{\|\mathbf{Y}^{(i)}\|}_F^2}}.
\end{equation}
Measuring the percentage of error of a low-rank approximation with respect to a large matrix has been popularly used in a variety of applications for dimensionality reduction such as Principle Component Analysis (PCA) \cite{Ref21}, ISOMAP learning \cite{Ref22}, Locally Linear Embedding (LLE) \cite{Ref23}, and Linear Discriminant Analysis (LDA) \cite{Ref24}. According to the Eckart-Young-Mirsky theorem \cite{Ref25}, for a fixed $k$ value, the optimal $k$-rank approximation has the minimal error percentage of $\mathbf{Y}^{(i)}$, which is
\begin{equation}
\epsilon_{min}={{{\|\mathbf{Y}_L^{(i)}-\mathbf{Y}^{(i)}\|}_F^2} \over {{\|\mathbf{Y}^{(i)}\|}_F^2}}={{\sum_{j=k+1}^{\min(m,n)}{\sigma_j^{(i)}}^2} \over {{\|\mathbf{Y}^{(i)}\|}_F^2}},
\end{equation}
where $\sigma_j$ is the $j$th singular value of $\mathbf{Y}^{(i)}$.

The adaptivity of R\textsuperscript{3}SVD is achieved by gradually constructing the low-rank approximation of $\mathbf{Y}^{(i)}$ while estimating the error percentage. The rationale of R\textsuperscript{3}SVD is to build a low-rank QB decomposition incrementally based on orthogonal Gaussian projection and then derive the low-rank SVD. There are a couple of attractive properties of QB decomposition:
\begin{enumerate}
\item[1).] Assuming that $\mathbf{QB}=\mathbf{Y}_L^{(i)}$ is a $k$-rank QB decomposition to approximate $\mathbf{Y}^{(i)}$, where $\mathbf{Q} \in \mathbb{R}^{m \times k}$ is orthonormal and $\mathbf{B} \in \mathbb{R}^{k \times n}$, then the error percentage can be evaluated as
\begin{equation}
{{{\|\mathbf{Y}_L^{(i)}-\mathbf{Y}^{(i)}\|}_F^2} \over {{\|\mathbf{Y}^{(i)}\|}_F^2}}={{{\|\mathbf{Y}^{(i)}\|}_F^2-{\|\mathbf{B}\|}_F^2} \over {{\|\mathbf{Y}^{(i)}\|}_F^2}}={{1-{\|\mathbf{B}\|}_F^2} \over {{\|\mathbf{Y}^{(i)}\|}_F^2}};
\end{equation}
\item[2).] Assuming that $\mathbf{Q}=[\mathbf{Q}_1,...,\mathbf{Q}_r]$ where $\mathbf{Q}_1,...,\mathbf{Q}_r$ are block row matrices and correspondingly, $\mathbf{B}=\left(\begin{matrix}\mathbf{B}_1 \\ \vdots \\ \mathbf{B}_r\end{matrix}\right)$, then, 
\begin{equation}
{{{\|\mathbf{Y}_L^{(i)}-\mathbf{Y}^{(i)}\|}_F^2} \over {{\|\mathbf{Y}^{(i)}\|}_F^2}}={{1-{\sum_j\|\mathbf{B}_j\|}_F^2} \over {{\|\mathbf{Y}^{(i)}\|}_F^2}}.
\end{equation}
\end{enumerate}
The mathematical proofs of the above properties of QB decomposition can be found in \cite{Ref6}. According to property 1), the error percentage of the low-rank approximation can be efficiently evaluated by computing the Frobenius norms of $\mathbf{B}$ and $\mathbf{Y}^{(i)}$. According to property 2), the error percentage can also be evaluated incrementally by adding up the Frobenius norm of $\mathbf{B}_j$ when incrementally building up the low-rank approximation.

By taking advantage of the nice properties of QB decomposition, we are able to design QB decomposition that can be built adaptively to satisfy a specific error percentage. Instead of using a fixed rank $k$, a rank value $r$ is adaptively derived. Initially, a $t$-rank QB decomposition is obtained and its error percentage is calculated accordingly, where $t$ is an initial guess of the appropriate rank $r$ and can also be justified according to the amount of memory available in the computer system. If the error percentage obtained so far does not satisfy the desired error percentage, a new $\Delta t$-rank QB approximation is built in the subspace orthogonal to the space of the previous QB approximation and is integrated with the previous QB decomposition. The error of the $\Delta t$-rank QB decomposition is also calculated and is then used to estimate the error percentage of the overall QB decomposition. The above process is repeated until the incrementally built low-rank approximation has error percentage less than the fixed precision threshold $\epsilon_{threshold}$. Finally, based on the obtained QB decomposition, a low-rank SVD decomposition satisfying the specific error percentage with its estimated rank $r$ is derived.

The matrix-block matrix multiplications between $\mathbf{Y}^{(i)}$ and the block matrices in R\textsuperscript{3}SVD take $O((np+2) r T_{mul})$ floating point operations. Here, $r$ is derived adaptively for a given precision. Assuming that, in a low-rank matrix where $m, n \gg r$, the computational cost of matrix-block matrix multiplications dominates those of QR, SVD, and orthogonalization operations, the main computational gain of R\textsuperscript{3}SVD compared to RSVD with an overestimating rank $k$ $(k>r)$ is the saving of $O((np+2)(k-r)T_{mul})$ floating operations.

\begin{algorithm}
 \KwIn{$\mathbf{Y}^{(i)} \in \mathbb{R}^{m \times n}$, initial sampling size $t \in \mathbb{N}$, sampling incremental step $\Delta t \in \mathbb{N}$ per iteration, number of power iteration $np \in \mathbb{N}$, and error percentage threshold $\epsilon_{threshod} \in \mathbb{R}$.}
 \KwOut{Low-rank approximation $\mathbf{U}_L^{(i)} \mathbf{\Sigma}_L^{(i)} {\mathbf{V}_L^{(i)}}^T$ with $\mathbf{U}_L^{(i)} \in \mathbb{R}^{m \times r}$, $\mathbf{\Sigma}_L^{(i)} \in \mathbb{R}^{r \times r}$, $\mathbf{V}_L^{(i)} \in \mathbb{R}^{r \times n}$, and estimated rank $r$}
 \bigbreak
 \tcc{build initial QB decomposition}
 Construct an $n \times t$ Gaussian random matrix $\mathbf{\Omega}$\;
 $\mathbf{X} \gets \mathbf{Y}^{(i)}\mathbf{\Omega}$\;
 \For {$j\leftarrow 1$ \KwTo $np$} { 
 $\mathbf{X} \gets \mathbf{Y}^{(i)}\mathbf{X}$ \tcc*{power iterations}
 }
 $[\mathbf{Q},\mathbf{R}] \gets qr(\mathbf{X})$ \tcc*{QR decomposition}
 $\mathbf{B} \gets \mathbf{Q}^T \mathbf{Y}^{(i)}$ \tcc*{QB decomposition}
 $normB \gets \|\mathbf{B}\|_F^2$\;
 $\epsilon \gets ({\|\mathbf{Y}^{(i)}\|}_F^2 - normB)/{\|\mathbf{Y}^{(i)}\|}_F^2$ \tcc*{error percentage}
 $r \gets t$\;
 \bigbreak
 \tcc{incrementally build up QB decomposition}
 \While {$\epsilon > \epsilon_{threshod}$} {
	Construct an $n \times \Delta t$ Gaussian random matrix $\mathbf{\Omega}$\;
	$\mathbf{X} \gets \mathbf{Y}^{(i)}\mathbf{\Omega}$\;
	\For {$j\leftarrow 1$ \KwTo $np$} { 
 	$\mathbf{X} \gets \mathbf{Y}^{(i)}\mathbf{X}$ \tcc*{power iterations}
 	}
	$\mathbf{X} \gets \mathbf{X}-\mathbf{QQ}^T \mathbf{X}$ \tcc*{orthogonalization with $\mathbf{Q}$}
	$[\mathbf{Q}',\mathbf{R}] \gets qr(\mathbf{X})$ \tcc*{QR decomposition}
 	$\mathbf{B}' \gets {\mathbf{Q}'}^T \mathbf{Y}^{(i)}$ \tcc*{QB decomposition}
 	$\mathbf{Q} \gets [\mathbf{Q}',\mathbf{Q}]$ \tcc*{build up approximate basis}
 	$\mathbf{B} \gets \left(\begin{matrix}\mathbf{B}' \\ \mathbf{B}\end{matrix}\right)$ \tcc*{gradually build up $\mathbf{B}$}
 	$normB \gets normB + \|\mathbf{B}'\|_F^2$\;
 	$\epsilon \gets ({\|\mathbf{Y}^{(i)}\|}_F^2 - normB)/{\|\mathbf{Y}^{(i)}\|}_F^2$ \tcc*{error percentage}
 	$r \gets r + \Delta t$\;
 }
 $[\mathbf{U}_L^{(i)}, \mathbf{\Sigma}_L^{(i)}, \mathbf{V}_L^{(i)}] \gets svd(\mathbf{B})$\;
 $\mathbf{U}_L^{(i)} \gets \mathbf{Q}\mathbf{U}_L^{(i)}$\;
 \caption{Rank Revealing Randomized SVD (R\textsuperscript{3}SVD)}
\end{algorithm}

\subsection{Singular Vectors Recycling}
\label{sec:2.3}
The SVT algorithm uses an iterative process to build up the low rank approximation, where SVD operation is repeatedly carried out on $\mathbf{Y}^{(i)}$. By observing the SVT process, one can find that the subsequent matrices $\mathbf{Y}^{(i)}$ and $\mathbf{Y}^{(i+1)}$ do not changes dramatically, particularly in the latter stage of the SVT process. Based on Equation (\ref{Bregman}), we have
\begin{equation}
\|\mathbf{Y}^{(i+1)}-\mathbf{Y}^{(i)}\|_F=\|\delta P_\mathbf{\Lambda}(\mathbf{A}-\mathbf{Y}^{(i)})\|_F \leq \delta\|\mathbf{A}-\mathbf{Y}^{(i)}\|_F.
\end{equation}
According to Theorem 4.1 in \cite{Ref1}, $\|\mathbf{Y}^{(i)} - \mathbf{A}\|_F \rightarrow 0$ as SVT is designed as a convex optimizer. Hence, the left singular vectors $\mathbf{U}^{(i)}$ of $\mathbf{Y}^{(i)}$ can be taken advantage as an approximate basis to capture a significant portion of actions of $\mathbf{Y}^{(i+1)}$ in the subsequent SVT iteration. As a consequence, the partial SVD operation on $\mathbf{Y}^{(i+1)}$ can be accelerated.

\subsection{Recycling Rank Revealing Randomized SVD (R\textsuperscript{4}SVD) Algorithm}
\label{sec:2.4}
At every SVT iteration except for the very first one, the R\textsuperscript{4}SVD algorithm reuses the left singular vectors of $\mathbf{Y}^{(i-1)}$ from the previous iteration. Compared to subspace reusing proposed by \cite{Ref7} in the FRSVT implementation, the reused subspace in R\textsuperscript{4}SVD does not involve in the power iterations, which is the main computational bottleneck of FRSVT \cite{Ref7}. Instead, R\textsuperscript{4}SVD uses the left singular vectors obtained from the previous iteration as the starting approximate basis and then incrementally builds up the QB decomposition satisfying the desired error percentage by constructing subspace orthogonal to the previous approximate basis. In the meanwhile, the approximation errors are monitored until the error percentage threshold is reached. Finally, carrying out SVD on the obtained QB decomposition provides the low-rank SVD approximation of $\mathbf{A}$ satisfying the specified error percentage threshold $\epsilon_{threshold}$. The complete R\textsuperscript{4}SVD algorithm is described in pseudocode as follows. 

\begin{algorithm}
 \KwIn{$\mathbf{Y}^{(i)} \in \mathbb{R}^{m \times n}$, left singular vectors from previous iteration $\mathbf{U}^{(i-1)} \in \mathbb{R}^{m \times s}$, sampling incremental step $\Delta t \in \mathbb{N}$, number of power iteration $np \in \mathbb{N}$, and error percentage threshold $\epsilon_{threshod} \in \mathbb{R}$.}
 \KwOut{Low-rank approximation $\mathbf{U}_L^{(i)} \mathbf{\Sigma}_L^{(i)} {\mathbf{V}_L^{(i)}}^T$ with $\mathbf{U}_L^{(i)} \in \mathbb{R}^{m \times r}$, $\mathbf{\Sigma}_L^{(i)} \in \mathbb{R}^{r \times r}$, $\mathbf{V}_L^{(i)} \in \mathbb{R}^{r \times n}$, and estimated rank $r$}
 \bigbreak
 \tcc{build initial QB decomposition}
 $\mathbf{Q} \gets \mathbf{U}^{(i-1)}$ \tcc*{Recycling}
 $\mathbf{B} \gets \mathbf{Q}^T \mathbf{Y}^{(i)}$ \tcc*{QB decomposition}
 $normB \gets \|\mathbf{B}\|_F^2$\;
 $\epsilon \gets ({\|\mathbf{Y}^{(i)}\|}_F^2 - normB)/{\|\mathbf{Y}^{(i)}\|}_F^2$ \tcc*{error percentage}
 $r \gets s$\;
 \bigbreak
 \tcc{incrementally build up QB decomposition}
 \While {$\epsilon > \epsilon_{threshod}$} {
	Construct an $n \times \Delta t$ Gaussian random matrix $\mathbf{\Omega}$\;
	$\mathbf{X} \gets \mathbf{Y}^{(i)}\mathbf{\Omega}$\;
	\For {$j\leftarrow 1$ \KwTo $np$} { 
 	$\mathbf{X} \gets \mathbf{Y}^{(i)}\mathbf{X}$ \tcc*{power iterations}
 	}
	$\mathbf{X} \gets \mathbf{X}-\mathbf{QQ}^T \mathbf{X}$ \tcc*{orthogonalization with $\mathbf{Q}$}
	$[\mathbf{Q}',\mathbf{R}] \gets qr(\mathbf{X})$ \tcc*{QR decomposition}
 	$\mathbf{B}' \gets {\mathbf{Q}'}^T \mathbf{Y}^{(i)}$ \tcc*{QB decomposition}
 	$\mathbf{Q} \gets [\mathbf{Q}',\mathbf{Q}]$ \tcc*{build up approximate basis}
 	$\mathbf{B} \gets \left(\begin{matrix}\mathbf{B}' \\ \mathbf{B}\end{matrix}\right)$ \tcc*{gradually build up $\mathbf{B}$}
 	$normB \gets normB + \|\mathbf{B}'\|_F^2$\;
 	$\epsilon \gets ({\|\mathbf{Y}^{(i)}\|}_F^2 - normB)/{\|\mathbf{Y}^{(i)}\|}_F^2$ \tcc*{error percentage}
 	$r \gets r + \Delta t$\;
 }
 $[\mathbf{U}_L^{(i)}, \mathbf{\Sigma}_L^{(i)}, \mathbf{V}_L^{(i)}] \gets svd(\mathbf{B})$\;
 $\mathbf{U}_L^{(i)} \gets \mathbf{Q}\mathbf{U}_L^{(i)}$\;
 \caption{Recycling Rank Revealing Randomized SVD (R\textsuperscript{4}SVD)}
\end{algorithm}

\subsection{Fast SVT Algorithm}
\label{sec:2.5}
At the $i$th SVT iteration, only the singular values greater than $\tau$ and their corresponding singular vectors are needed in $D_{\tau}(\mathbf{Y}^{(i)})$. The number of singular values needed increases gradually as SVT iterates. At the early stage of estimating $D_{\tau}(\mathbf{Y}^{(i)})$, there are often only a few singular values exceeding $\tau$. Therefore, the QB decomposition with high error percentage is sufficient to estimate these singular values and their associated singular vectors with desired accuracy. Nevertheless, as SVT progresses, the rank of the approximated low-rank decomposition increases gradually, which demands the error percentage of QB decomposition to decrease accordingly. 

Here, we employ a simulated annealing style cooling scheme in the fast SVT implementation to adaptively adjust the error percentage threshold $\epsilon_{threshod}$. We keep track of the approximation error $\epsilon^{(i)}$ defined as
\begin{equation}
\epsilon^{(i)}=\|P_\mathbf{\Lambda}(\mathbf{A}-\mathbf{Y}^{(i)})\|.
\end{equation}
Due to the fact that SVT is a convex optimization algorithm, $\epsilon^{(i)}$ is supposed to decrease continuously. Once $\epsilon^{(i)}$ stops decreasing, it indicates that the error in the approximated partial SVD operation is too high. This triggers the annealing operation to reduce the error percentage threshold parameter $\epsilon_{threshold}$ for the subsequent R\textsuperscript{4}SVD operations such that 
\begin{equation}
\epsilon_{threshold}=\beta \epsilon_{threshold},
\end{equation}
where $0<\beta<1$ is the annealing factor. 

By putting all puzzles together, the fast implementation of SVT algorithm using R\textsuperscript{4}SVD to perform partial SVD operation is described as follows.

\begin{algorithm}
 \KwIn{Sample set $\mathbf{\Lambda}$ and sampled entries $P_\mathbf{\Lambda}(\mathbf{A})$, step size $\delta \in \mathbb{R}$, maximum number of iterations $maxit \in \mathbb{N}$, sampling incremental step $\Delta t \in \mathbb{N}$, number of power iteration $np \in \mathbb{N}$, error tolerance $\epsilon \in \mathbb{R}$, and threshold parameter $\tau \in \mathbb{R}$.}
 \KwOut{Low-rank approximation $\mathbf{U}_L^{(i)} \mathbf{\Sigma}_L^{(i)} {\mathbf{V}_L^{(i)}}^T$ with $\mathbf{U}_L^{(i)} \in \mathbb{R}^{m \times r}$, $\mathbf{\Sigma}_L^{(i)} \in \mathbb{R}^{r \times r}$, $\mathbf{V}_L^{(i)} \in \mathbb{R}^{r \times n}$, and estimated rank $r$}
 \bigbreak
 \tcc{initialization}
 $\mathbf{Y}^{(0)} \gets \lceil {\tau \over {\delta\|P_\mathbf{\Lambda}(\mathbf{A})\|}} \rceil \delta \|P_\mathbf{\Lambda}(\mathbf{A})\|$\;
 $\mathbf{Q} \gets \varnothing, \mathbf{B} \gets \varnothing$\;
 $\mathbf{U}^{(0)}_L \gets \varnothing$\;
 $\epsilon_{min} \gets \infty$\;
 
 \bigbreak
 \tcc{SVT iterations}
 \For {$i\leftarrow 1$ \KwTo $maxit$} {
 	$[\mathbf{U}_L^{(i)}, \mathbf{\Sigma}_L^{(i)}, \mathbf{V}_L^{(i)}] \gets$ R\textsuperscript{4}SVD($\mathbf{Y}^{(i)}$, $\mathbf{U}_L^{(i-1)}$, $\Delta t$, $np$, $\epsilon_{threshold}$)\tcc*{R\textsuperscript{4}SVD}
	$\epsilon^{(i)} \gets \|P_\mathbf{\Lambda}(\mathbf{A}-\mathbf{Y}^{(i)})\|$\;
	\eIf {$\epsilon^{(i)} < \epsilon_{min}$} {
		$\epsilon_{min} \gets \epsilon^{(i)}$\;
	} {
		$\epsilon_{threshold} \gets \beta \epsilon_{threshold}$\tcc*{simulated annealing style cooling}
	}
	\If {$\epsilon^{(i)} < \epsilon$} {
		return $\mathbf{U}_L^{(i)} \mathbf{\Sigma}_L^{(i)} {\mathbf{V}_L^{(i)}}^T$\;
	}
	$\mathbf{X}^{(i)} \gets D_{\tau}(\mathbf{U}_L^{(i)}, \mathbf{\Sigma}_L^{(i)}, \mathbf{V}_L^{(i)})$\tcc*{linearized Bregman’s iterations}
	$\mathbf{Y}^{(i+1)} \gets \mathbf{Y}^{(i)} + \delta P_\mathbf{\Lambda}(\mathbf{A} - \mathbf{X}^{(i)})$\;
 }
 \caption{Fast Singular Value Thresholding}
\end{algorithm}

\section{Numerical Results}
\label{sec:3}
We compare the performance of several SVT implementations using SVT-R\textsuperscript{4}SVD, SVT-Lanczos \cite{Ref1}, SVT-full SVD \cite{Ref1}, and FRSVT \cite{Ref7} on two matrix completion applications, including image recovery and movie recommendation. The initial sampling size $t$ and the sampling incremental step $\Delta t$ in SVT-R\textsuperscript{4}SVD are set to $\lfloor 0.05 \min (m,n) \rfloor$ and $10$, respectively. The annealing factor $\beta$ in SVT-R\textsuperscript{4}SVD is $0.95$. SVT-Lanczos is an SVT implementation based on PROPACK \cite{Ref11} using the Lanczos procedure with partial re-orthogonalization. SVT-full SVD carries out full SVD in every SVT iteration step. FRSVT is implemented according to Algorithm 1 in \cite{Ref7}. In all SVT implementations, threshold parameter $\tau$ is set to $\|P_\mathbf{\Lambda}(\mathbf{A})\|_F$ and step size $\delta$ is $\sqrt{m \times n / ns}$ as suggested by \cite{Ref1}, where $ns$ is the total number of known entries in $\mathbf{\Lambda}$. The computations are carried out on a Linux server with 40 2.3GHz, 10-core Xeon processors and 512GB memory using Matlab version 2016a. The measured CPU time is obtained by the Matlab cputime() function, which is irrespective of the number of threads and cores used in different SVT implementations.

\subsection{Image Recovery}
\label{sec:3.1}
\begin{figure}
  \captionsetup[subfigure]{justification=centering}
  \centering
  \subfloat[$512 \times 512$ Original]{\includegraphics[width=0.49\textwidth]{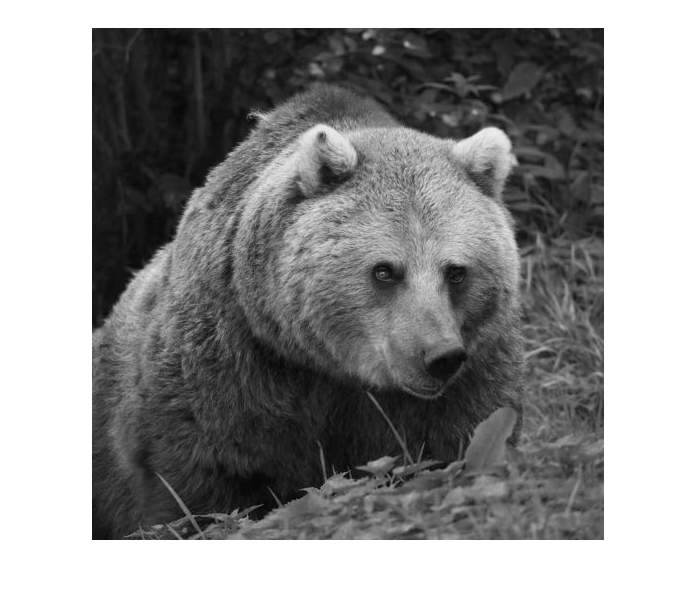}\label{fig:1a}}
  \hfill
  \subfloat[20\% Uniform Pixel Samples from Original Image]{\includegraphics[width=0.49\textwidth]{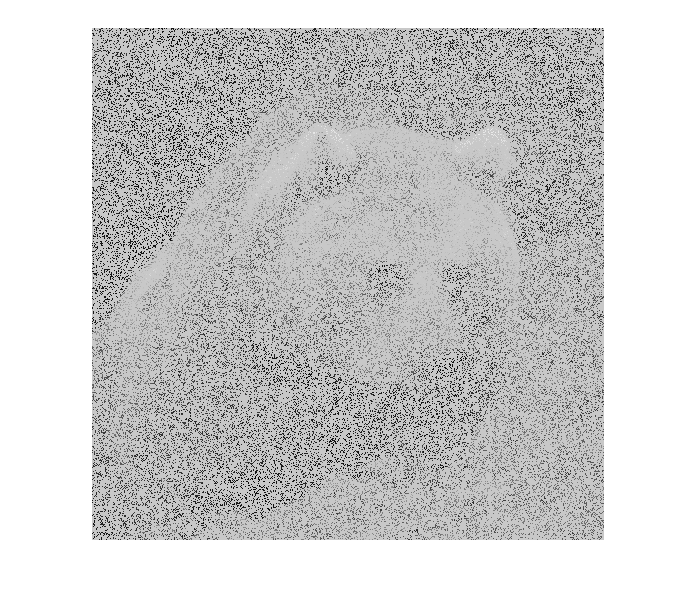}\label{fig:1b}}
  \hfill
  \subfloat[SVT-R\textsuperscript{4}SVD \newline(Rank: $142$, MAE: $17.22$, CPUTime: $211.6$s)
]{\includegraphics[width=0.49\textwidth]{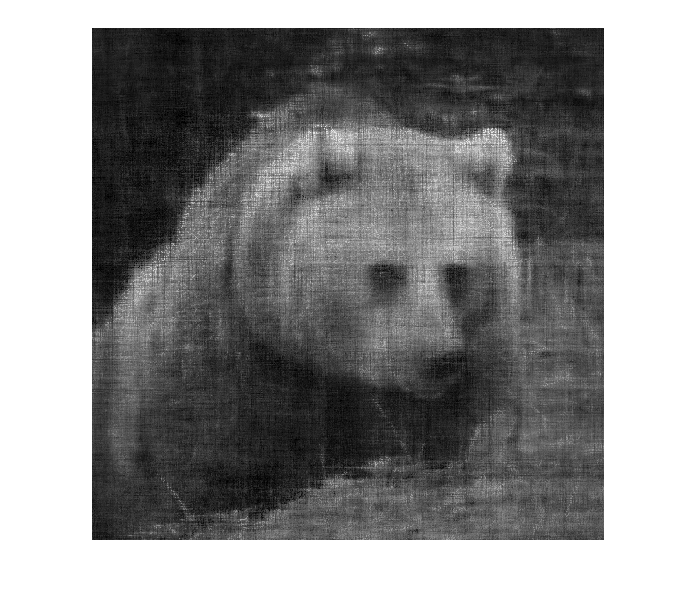}\label{fig:1c}}  
  \hfill
  \subfloat[SVT-Lanczos \newline(Rank: $141$, MAE: $17.19$, CPUTime: $1,038.1$s)]{\includegraphics[width=0.49\textwidth]{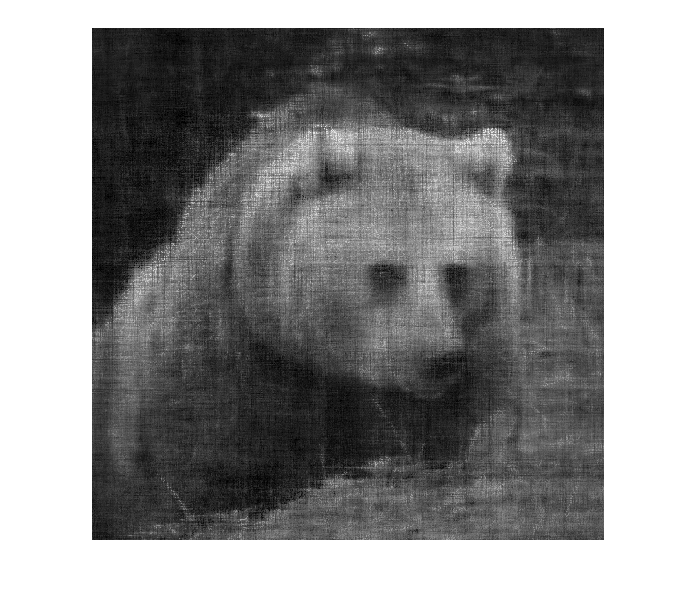}\label{fig:1d}}
  \hfill
  \subfloat[SVT-Full SVD \newline(Rank: $141$, MAE: $17.19$, CPUTime: $445.9$s)]{\includegraphics[width=0.49\textwidth]{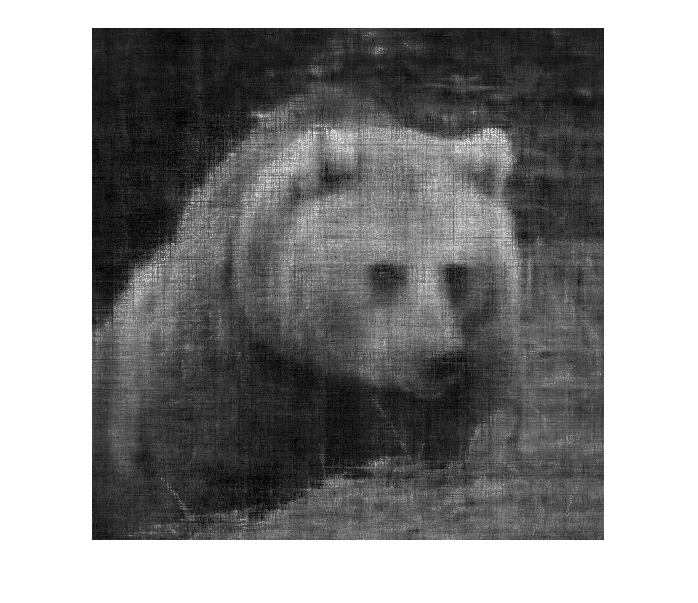}\label{fig:1e}}
  \hfill
  \subfloat[FRSVT \newline(Rank: $141$, MAE: $17.19$, CPUTime: $316.5$s)]{\includegraphics[width=0.49\textwidth]{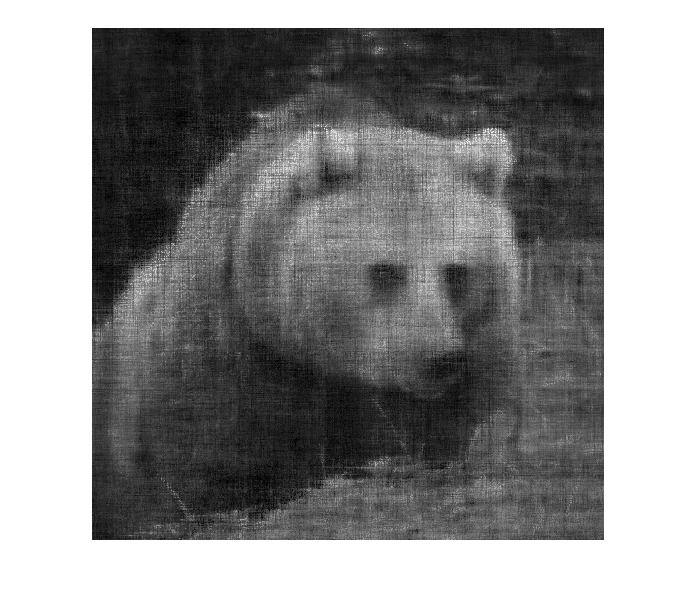}\label{fig:1f}}    
  \caption{Image recovery using SVT-R\textsuperscript{4}SVD, SVT-Lanczos, SVT-full SVD, and FRSVT on 20\% pixel samples of a $512 \times 512$ image.}
\label{fig:1}
\end{figure}

Fig. \ref{fig:1} compares the performance of SVT-R\textsuperscript{4}SVD, SVT-Lanczos, SVT-full SVD, and FRSVT in recovering a small $512 \times 512$ image (Fig. \ref{fig:1}(\subref*{fig:1a})) from 20\% uniformly distributed pixel samples (Fig. \ref{fig:1}(\subref*{fig:1b})). The image is represented as a $512 \times 512$ matrix whose entries are the $8$-bit grayscale values ranging from $0$ to $255$. The overall image recovery error is then measured by the Mean Absolute Error (MAE) defined as
\begin{equation}
{1 \over mn} \sum_{(i,j)}|\mathbf{A}_{ij}-{\mathbf{Y}_L}_{ij}|,
\end{equation}
where $\mathbf{Y}_L$ is the obtained low-rank approximation obtained by SVT. The SVT termination condition is set to when the MAE on the sample set is less than $1.0$, i.e, ${1 \over ns} \sum_{(i,j) \in \mathbf{\Lambda}}|\mathbf{A}_{ij}-{\mathbf{Y}_L}_{ij}| < 1.0$. The images recovered from various SVT implementations are shown in Fig. \ref{fig:1}(\subref*{fig:1c})-\ref{fig:1}(\subref*{fig:1f}). The image recovered by SVT-R\textsuperscript{4}SVD yields slightly larger recovery error (0.17\% bigger error compared to that of SVT-full SVD) and slightly bigger matrix rank than those from the other SVT implementations, due to the approximation error in partial SVD. However, the CPU time of SVT-R\textsuperscript{4}SVT is significantly reduced, which is less than half of the CPU time of SVT-full SVD and about $2/3$ of that of FRSVT. The CPU time of SVT-Lanczos is even much higher than that of the SVT implementation with full SVD, indicating that SVT-Lanczos is not suitable for completing small matrices.

\begin{figure}
  \centering
  \includegraphics[width=1.00\textwidth]{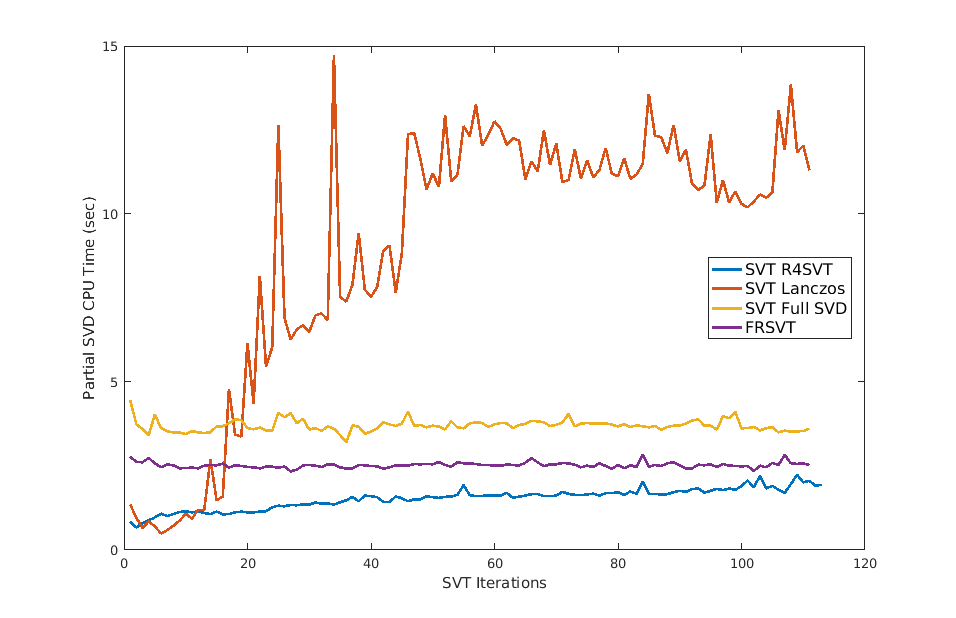}
  \caption{Comparison of CPU times of (partial) SVD operations in SVT-R\textsuperscript{4}SVD, SVT-Lanczos, SVT-full SVD, and FRSVT at every SVT iteration step when recovering a $512 \times 512$ image.}
\label{fig:2}
\end{figure}

Fig. \ref{fig:2} provides a comparison of CPU times of (partial) SVD operations in SVT-R\textsuperscript{4}SVD, SVT-Lanczos, SVT-full SVD, and FRSVT at every SVT iteration step. Due to the fact that full SVD is carried out on full rank and FRSVT is based on a fixed-rank SVD approximation, the CPU times of SVD operations at each SVT iteration step remain almost constant in SVT-full SVD and FRSVT. Partial SVD approximation using Lanczos algorithm in SVT-Lanczos is fast at the very beginning, but its computational cost increases rapidly after $20$ steps, even significantly higher than that of computing the full SVD, which agrees with the observation found in \cite{Ref5}. In comparison, SVT-R\textsuperscript{4}SVD is an adaptive partial SVD approximation algorithm where the size of the sampling matrix is adjusted adaptively according to the number of singular values/vectors needed to approximate the SVT operator at every step. As a result, the CPU time of partial SVD operations increases gradually but remains lower than that of FRSVT, due to avoidance of costly power iterations.

\begin{figure}
  \captionsetup[subfigure]{justification=centering}
  \centering
  \subfloat[$8,192 \times 8,192$ Original]{\includegraphics[width=0.435\textwidth]{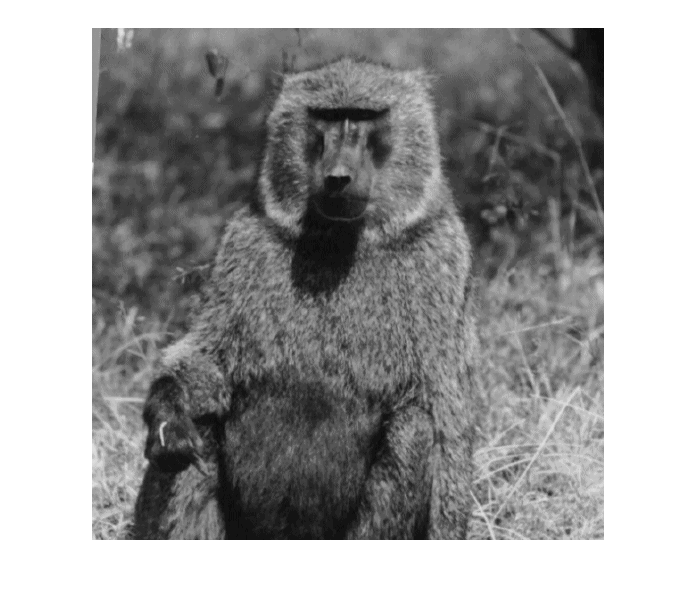}\label{fig:3a}}
  \hfill
  \subfloat[20\% Uniform Pixel Samples from Original Image]{\includegraphics[width=0.49\textwidth]{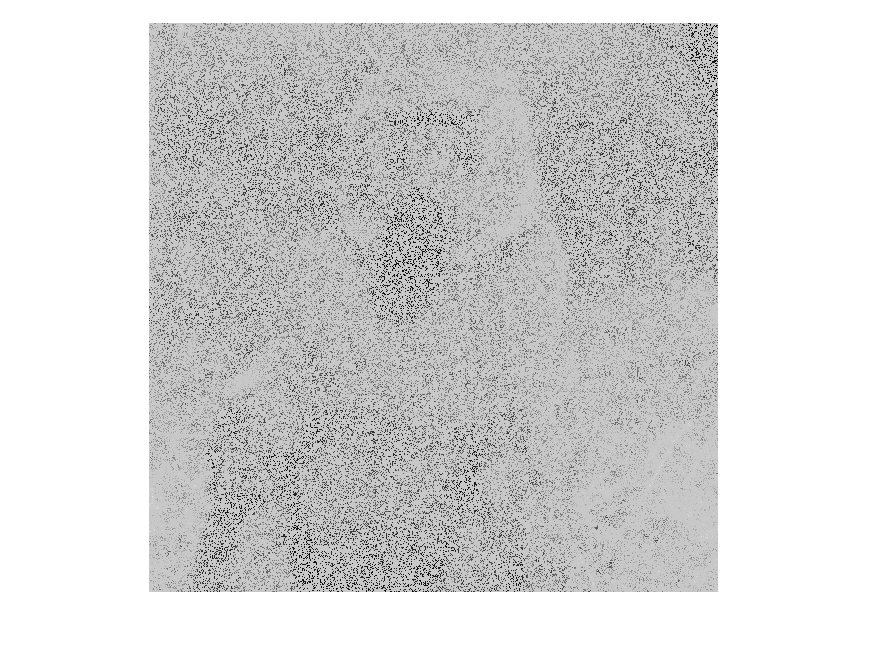}\label{fig:3b}}
  \hfill
  \subfloat[SVT-R\textsuperscript{4}SVD\newline(Rank: $364$, MAE: $2.34$, CPUTime: $16,223.0$s)
]{\includegraphics[width=0.49\textwidth]{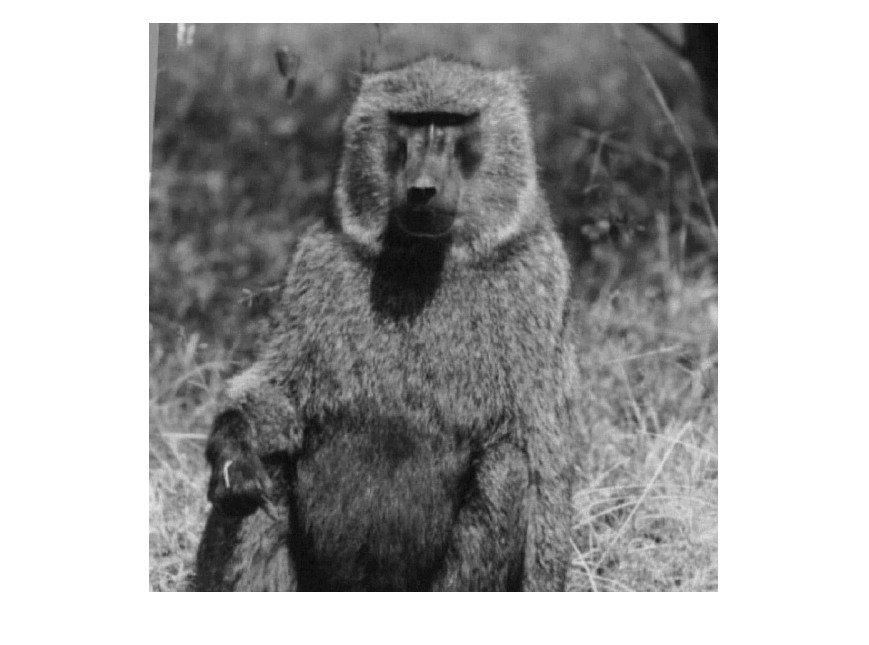}\label{fig:3c}}  
  \hfill
  \subfloat[SVT-Lanczos\newline(Rank: $341$, MAE: $2.32$, CPUTime: $87,496.8$s)]{\includegraphics[width=0.49\textwidth]{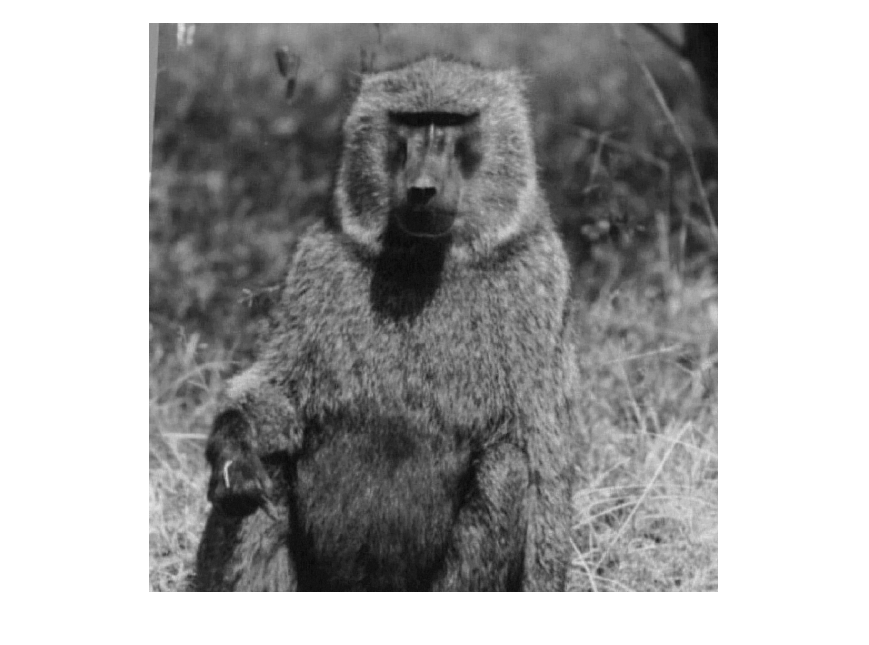}\label{fig:3d}}
  \hfill
  \subfloat[SVT-Full SVD\newline(Rank: $341$, MAE: $2.32$, CPUTime: $234,963.9$s)]{\includegraphics[width=0.49\textwidth]{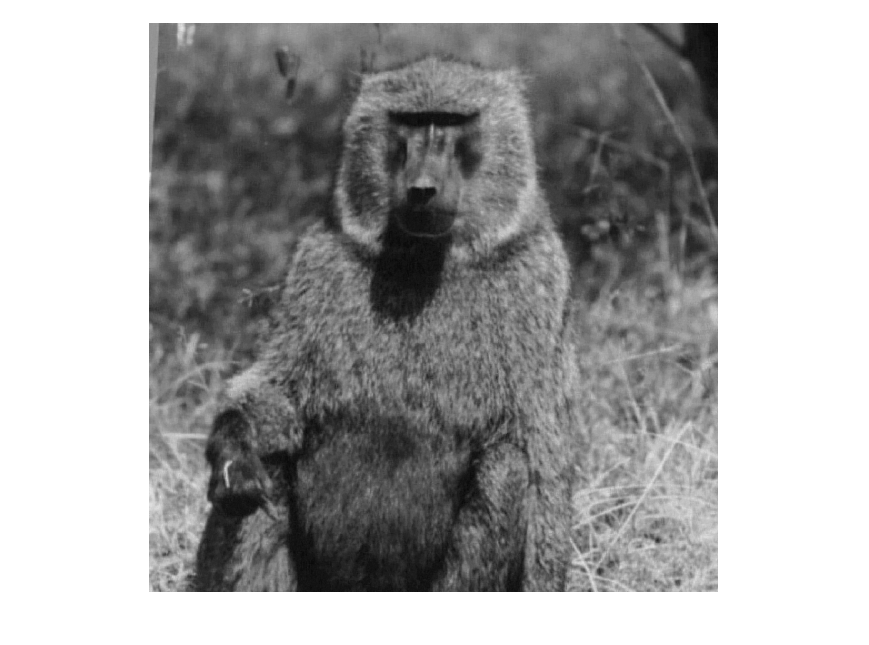}\label{fig:3e}}
  \hfill
  \subfloat[FRSVT\newline(Rank: $342$, MAE: $2.33$, CPUTime: $42,120.9$s)]{\includegraphics[width=0.49\textwidth]{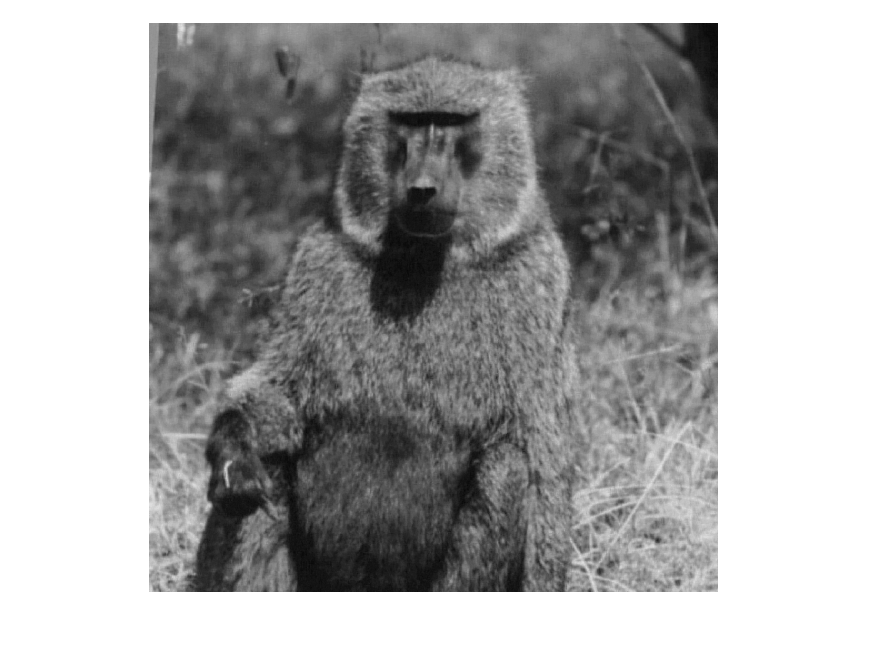}\label{fig:3f}}    
  \caption{Image recovery using SVT-R\textsuperscript{4}SVD, SVT-Lanczos, SVT-full SVD, and FRSVT on 20\% pixel samples of an $8,192 \times 8,192$ image.}
\label{fig:3}
\end{figure}

\begin{figure}
  \centering
  \includegraphics[width=1.00\textwidth]{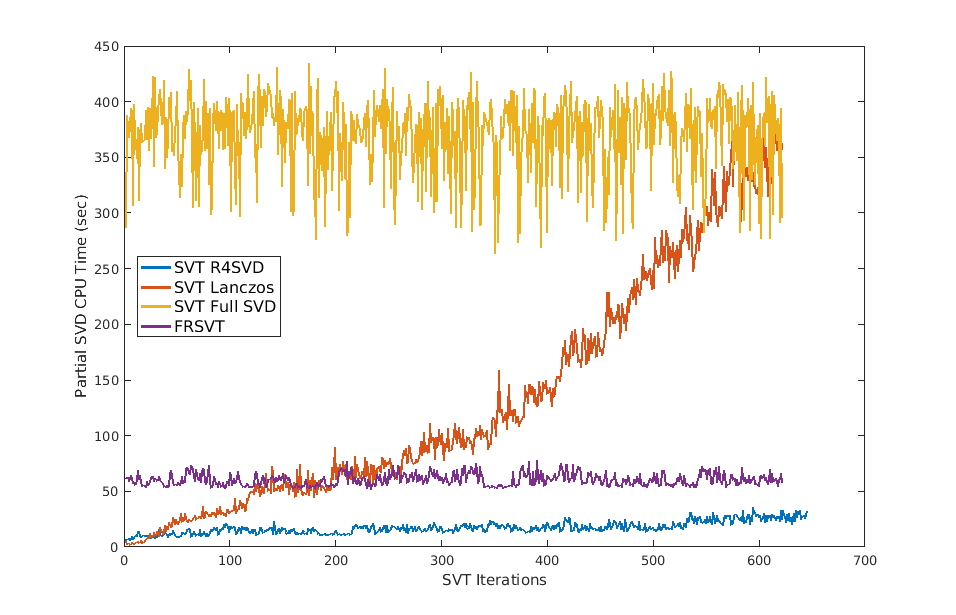}
  \caption{Comparison of CPU times of (partial) SVD operations in SVT-R\textsuperscript{4}SVD, SVT-Lanczos, SVT-full SVD, and FRSVT at every SVT iteration step when recovering an $8,192 \times 8,192$ image.}
\label{fig:4}
\end{figure}

Fig. \ref{fig:3} compares the performance of SVT-R\textsuperscript{4}SVD, SVT-Lanczos, SVT-full SVD, and FRSVT in recovering a big $8,192 \times 8,192$ image Fig. \ref{fig:3}(\subref*{fig:3a}) with 20\% uniformly distributed samples showed in Fig. \ref{fig:3}(\subref*{fig:3b}). Similar to recovering the small one, the image recovered by SVT-R\textsuperscript{4}SVD (Fig. \ref{fig:3}(\subref*{fig:3c})) yields slightly larger recovery error (0.86\% more error compared to that of SVT-full SVD) and slightly bigger recovery matrix rank than those by SVT-full SVD, SVT-Lanczos, and FRSVT (Figs. \ref{fig:3}(\subref*{fig:3d})-\ref{fig:3}(\subref*{fig:3f})). However, the CPU time of SVT-R\textsuperscript{4}SVD is only 6.9\%, 18.5\%, and 38.5\% of those of SVT-full SVD, SVT-Lanczos, and FRSVT, respectively. Fig. \ref{fig:4} compares the CPU times at every SVT iteration step. For a large matrix, full SVD becomes very costly. Partial SVD based on Lanczos algorithm is fast initially but its computational cost eventually catches up with that of full SVD in later SVT iterations. FRSVT and SVT-R\textsuperscript{4}SVD using partial SVD based on randomized sampling demonstrate clear advantages. Compared to FRSVT, SVT-R\textsuperscript{4}SVD yields more significant computational gains due to two reasons. One is that R\textsuperscript{4}SVD uses adaptive sampling, which leads to computational saving at the early stage of SVT when the number of singular values/vectors needed to approximate the SVT operator is small. The other reason is that R\textsuperscript{4}SVD avoids the costly power iterations in FRSVT.

It is important to notice that FRSVT uses a fixed-rank randomized SVD approximation, which requires good estimation of the appropriate rank of the recovered matrix. In the above computational results showed in recovering the $512 \times 512$ and $8,192 \times 8,192$ images, FRSVT uses fixed ranks of $150$ and $400$, respectively, which are close to the final recovered matrix rank, due to the fact that we know the appropriate ranks. However, in practice, such good estimations of the final matrix rank are often difficult to obtain beforehand. Consequently, an underestimation of the fixed rank will lead to SVT divergence while an overestimation will result in unnecessary computations. In contrast, R\textsuperscript{4}SVD is a fixed precision RSVD algorithm while SVT-R\textsuperscript{4}SVD adaptively adjusts partial SVD precision requirements along SVT iterations, where estimating a fixed recovery matrix rank beforehand is not necessary.

It is also interesting to notice that with 20\% pixel samples, the recovery error on the large image is much smaller than that on the small one. This is due to the fact that 20\% samples in the large image has 256 times more pixels than the small one while the actual rank ratio between the two images is less than 3.0. Therefore, with significantly more information, the SVT algorithm recovers the large image with better accuracy.

\subsection{Movie Recommendation}
\label{sec:3.2}
Here we apply various SVT implementations to a movie recommendation problem. Considering an incomplete matrix with rows as users, columns as movies, and entries as ratings, the goal of a movie recommendation system is to predict the blanks in the user-movie matrix to infer the unknown preferences of users on movies. We use a large user-movie matrix from MovieLens \cite{Ref12} containing $10,000,054$ ratings ($1$-$5$) from $71,567$ users on $65,133$ movie tags. We randomly select $80\%$ of the ratings to construct a training set and use the rest $20\%$ as the testing set. We use various SVT implementations to complete the user-movie matrix constructed from the training set and then validate the results using the testing set. The SVT termination condition is set to when MAE on the training set is less than $0.1$. Since carrying out full SVD on this large matrix is very costly, we only compare the performance of SVT-R\textsuperscript{4}SVD, SVT-Lanczos, and FRSVT here. 
	
\begin{figure}
  \centering
  \includegraphics[width=1.00\textwidth]{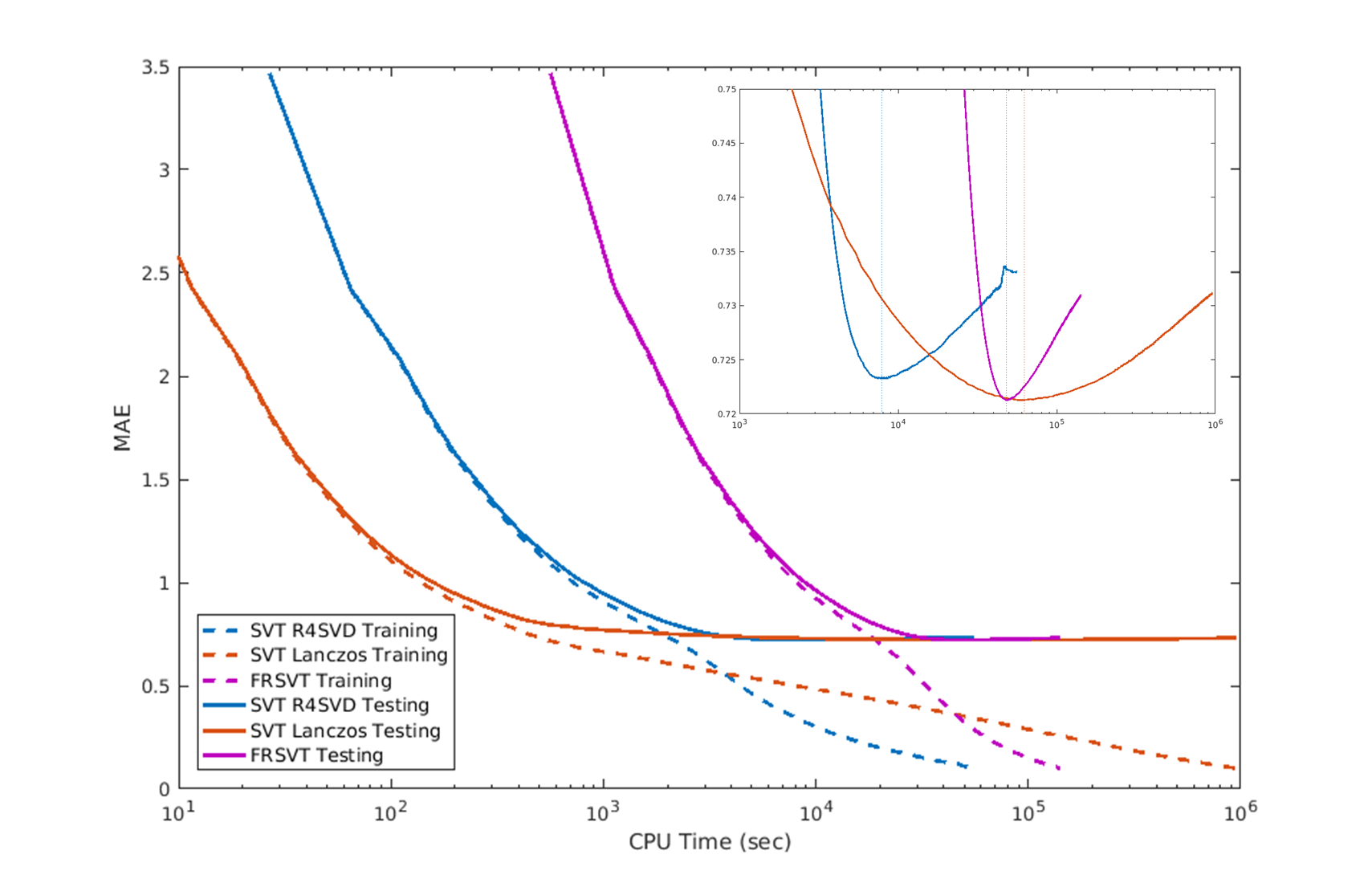}
  \caption{Performance comparison of movie recommendation using SVT-R\textsuperscript{4}SVD, SVT-Lanczos, and FRSVT on the $71,567 \times 65,133$ user-movie matrix from MovieLens.}
\label{fig:5}
\end{figure}

Fig. \ref{fig:5} plots the MAEs on the training set as well as the testing set with respect to the CPU times spent by SVT-R\textsuperscript{4}SVD, SVT-Lanczos, and FRSVT. For the training set, SVT-R\textsuperscript{4}SVD takes $55,878$ seconds of CPU time to reach MAE $< 0.1$, which is $39.4\%$ of that of FRSVT ($141,624$s) and $5.8\%$ of that of SVT-Lanczos ($956,746$s). The MAEs on the testing set cannot reach $0.1$. Instead, overfitting starts to occur after certain SVT iteration steps - the testing set MAE curves reaches a minimum and then slowly increases, as shown in the inner figure in Fig. \ref{fig:5}. After all, SVT-R\textsuperscript{4}SVD reaches the optimal testing results at $7,867$s while in comparison, FRSVT and SVT-Lanczos obtains the optimal testing results at $49,432$s and $62,425$s. The tradeoff is, SVT-R\textsuperscript{4}SVD yields slightly higher optimal testing MAE ($0.7228$) compared to $0.7214$ in FRSVT and $0.7212$ in SVT-Lanczos. 

\section{Conclusions}
\label{sec:4}
In this paper, a fast implementation of the SVT algorithm is developed for matrix completion applications. We propose an SVT-R\textsuperscript{3}SVD algorithm to fast approximate the SVT operator so as to satisfy a given fixed precision. We also show that the singular vectors obtained from the previous SVT iteration can be recycled as the approximate basis in the subsequent iteration to reduce the computational cost of partial SVD when estimating the SVT operator. The corresponding algorithms is so-called the SVT-R\textsuperscript{4}SVD algorithm. The SVT-R\textsuperscript{4}SVD algorithm adopts a simulated annealing style cooling mechanism to adaptively adjust the low-rank approximation error threshold along SVT iterations. The effectiveness of SVT-R\textsuperscript{4}SVD, in comparison with SVT-Lanczos, SVT-full SVD, and FRSVT, is demonstrated in matrix completion applications including image recovery and movie recommendation systems. The adaptiveness of SVT-R\textsuperscript{4}SVD leads to significant computation savings in completing both small and large matrices, with very small scarification of accuracy. 

\begin{acknowledgements}
Yaohang Li acknowledges support from National Science Foundation Grant Number 1066471.
\end{acknowledgements}


\begin{thebibliography}{}
\bibitem{Ref1}
J. Cai, E. J. Candes, Z. Shen, “A singular value thresholding algorithm for matrix completion,” SIAM J. Optim., 20(4): 1956-1982, 2010.
\bibitem{Ref2}
E. Candes, B. Recht, “Simple bounds for recovering low-complexity models,” Mathematical Programming, 141(1-2): 577-589, 2013.
\bibitem{Ref3}
J. Cai, S. Osher, “Fast singular value thresholding without singular value decomposition,” Methods and Applications of Analysis, 20(4): 335-352, 2013.
\bibitem{Ref4}
R. Larsen. PROPACK: Software for large and sparse SVD calculations. http://soi.stanford.edu/rmunk/PROPACK, 2017.
\bibitem{Ref5}
Z. Lin, M. Chen, Y. Ma, “The augmented Lagrange multiplier method for exact recovery of corrupted low-rank matrices,” arXiv:1009.5055, 2010.
\bibitem{Ref6}
Y. Gu, W. Yu, Y. Li, “Efficient randomized algorithms for adaptive low-rank factorizations of large matrices,” arXiv: 1606.09402, 2016.
\bibitem{Ref7}
T. Oh, Y. Matsushita, Y. Tai, I. Kweon, “Fast Randomized Singular Value Thresholding for Nuclear Norm Minimization,” Proceedings of CVPR2015, 2015.
\bibitem{Ref8}
N. Halko, P. G. Martinsson, J. A. Tropp, “Finding Structure with Randomness: Probabilistic Algorithms for Constructing Approximate Matrix Decompositions,” SIAM Rev., 53(2): 217–288, 2009.
\bibitem{Ref9}
K. J. Arrow, L. Hurwicz, H. Uzawa, Studies in Linear and Nonlinear Programming, Stanford University Press, Stanford, CA, 1958.
\bibitem{Ref10}
W. Yin, S. Osher, J. Darbon, D. Goldfarb. "Bregman Iterative Algorithms for Com-pressed Sensing and Related Problems." SIAM Journal on Imaging Sciences, 1(1):143-168, 2008.
\bibitem{Ref11}
R. Larsen, PROPACK: Software for large and sparse SVD calculations. http://soi.stanford.edu/rmunk/PROPACK.
\bibitem{Ref12}
F. M. Harper, J. A. Konstan, “The MovieLens Datasets: History and Context,” ACM Transactions on Interactive Intelligent Systems, 5(4): 19, 2015.
\bibitem{Ref13}
S. Kirkpatrick, C.D. Gelatt Jr., M.P. Vecchi, “Optimization by simulated annealing,” Science, 220: 671–680, 1983.
\bibitem{Ref14}
Y. Li, V. A. Protopopescu, N. Arnold, X. Zhang, A. Gorin, “Hybrid Parallel Temper-ing/Simulated Annealing Method,” Applied Mathematics and Computation, 212: 216-228, 2009.
\bibitem{Ref15}
S. Ma, D. Goldfarb, L. Chen. “Fixed point and Bregman iterative methods for matrix rank minimization,” Mathematical Programming, 128(1):321-353, 2011.
\bibitem{Ref16}
Y. Mu, J. Dong, X. Yuan, S. Yan. “Accelerated low-rank visual recovery by random projection,” Proceedings of IEEE Conference on Computer Vision and Pattern Recognition (CVPR), 2011.
\bibitem{Ref17}
S. Boyd, L. Vandenberghe. “Convex optimization,” Cambridge University press, 2004. 
\bibitem{Ref18}
Candes, X. Li, Y. Ma, J. Wright, “Robust principal component analysis,” Journal of the ACM, 58(3):11, 2011.
\bibitem{Ref19}
Y. Liu, D. Sun, K. Toh, “An implementable proximal point algorithmic framework for nuclear norm minimization,” Mathematical Programming, 133(1):399–436, 2012.
\bibitem{Ref20}
K. Toh, S. Yun, “An accelerated proximal gradient algorithm for nuclear norm regular-ized least squares 
problems,” Pacific J. Optimization, 6(3), 2010.
\bibitem{Ref21}
I. Jolliffe, Principal Component Analysis, 2nd ed. New York, NY, USA:Springer-Verlag, 2002. 
\bibitem{Ref22}
J. B. Tenenbaum, V. de Silva, and J. C. Langford, “A global geometric framework for nonlinear dimensionality reduction,” Science, vol. 290, no. 5500, pp. 2319–2323, Dec. 2000
\bibitem{Ref23}
S. T. Roweis and L. K. Saul, “Nonlinear dimensionality reduction by locally linear embedding,” Science, vol. 290, no. 5500, pp. 2323–2326, Dec. 2000.
\bibitem{Ref24}
S. Mika, G. Rätsch, J. Weston, B. Schölkopf, and K.-R. Müller, “Fisher discriminant analysis with kernels,” in Proc. Neural Networks for Signal Processing IX, Y.-H. Hu, J. Larsen, E. Wilson, and S. Douglas, Eds., pp. 41–48, 1999.
\bibitem{Ref25}
C. Eckart and G. Young, “The approximation of one matrix by another of lower rank,” Psychometrica, vol. 1, pp. 211–218, 1936.
\bibitem{Ref26}
A. K. Menon, and C. Elkan, “Fast Algorithms for Approximating the Singular Value Decomposition,” ACM Transactions on Knowledge Discovery from Data, vol. 5, No. 2, pp. 13:1-13:36, 2011.  

% etc
\end{thebibliography}


\end{document}